\documentclass[11pt]{article}
\usepackage[height=8in]{geometry}
\usepackage{amsmath, amssymb}
\usepackage{graphicx, color}
\usepackage{array, multirow}

\usepackage[sort]{natbib}

\setcitestyle{numbers,square,comma}
\usepackage[font=small,labelfont=bf]{caption} 

\usepackage{hyperref}
\hypersetup{colorlinks=true,
			urlcolor=blue,
			linkcolor=blue,
			citecolor=blue,
			bookmarksdepth=paragraph,
			pdftitle={Bounding extrema over global attractors using polynomial optimization (arXiv version)},
			pdfauthor={David Goluskin},}
\usepackage[nameinlink]{cleveref}
\crefname{equation}{}{}
\crefname{figure}{figure}{figures}
\crefname{table}{table}{tables}
\crefname{section}{section}{sections}
\Crefname{section}{Section}{Sections}

\definecolor{mygray}{rgb}{.6, .6, .6}

\newcommand\solidrule[1][15pt]{\rule[0.5ex]{#1}{1pt}}
\newcommand\dashedrule{\mbox{%
	\solidrule[3pt]\hspace{3pt}\solidrule[3pt]\hspace{3pt}\solidrule[3pt]}}
\newcommand\longdashedrule{\mbox{%
	\solidrule[6pt]\hspace{2pt}\solidrule[6pt]}}
\newcommand\dottedrule{\mbox{%
	\solidrule[1.5pt]\hspace{2pt}\solidrule[1.5pt]\hspace{2pt}\solidrule[1.5pt]\hspace{2pt}\solidrule[1.5pt]\hspace{2pt}\solidrule[1.5pt]}}
\newcommand\dashdottedrule{\mbox{%
	\solidrule[4.5pt]\hspace{2pt}\solidrule[2pt]\hspace{2pt}\solidrule[4.5pt]}}

\newcommand{\beq}{\begin{equation}}
\newcommand{\eeq}{\end{equation}}

\newcommand{\x}{\mathbf{x}}

\newcommand{\f}{\mathbf{f}}

\newcommand{\ddt}{\frac{\rm d}{{\rm d}t}}
\newcommand{\R}{\mathbb{R}}

\newcommand{\cA}{\mathcal{A}}
\newcommand{\eps}{\varepsilon}
\newcommand{\Rey}{{\rm Re}}

\title{Bounding extrema over global attractors\\using polynomial optimization}
\author{David Goluskin\thanks{Email: {\tt goluskin@uvic.ca}}}
\date{Department of Mathematics and Statistics, University of Victoria, Canada}

\begin{document}

\maketitle

\begin{abstract}
We present a framework for bounding extreme values of quantities on global attractors of differential dynamical systems. A global attractor is the minimal set that attracts all bounded sets; it contains all forward-time limit points. Our approach uses (generalized) Lyapunov functions to find attracting sets, which must contain the global attractor, and the choice of Lyapunov function is optimized based on the quantity whose extreme value one aims to bound. We also present a non-global framework for bounding extrema over the minimal set that is attracting in a specified region of state space. If the dynamics are governed by ordinary differential equations, and the equations and quantities of interest are polynomial, then our methods can be implemented computationally using polynomial optimization. In particular, we enforce nonnegativity of certain polynomial expressions by requiring them to be representable as sums of squares, leading to a convex optimization problem that can be recast as a semidefinite program and solved computationally. This computer assistance lets one construct complicated polynomial Lyapunov functions. Computations are illustrated using three examples. The first is the chaotic Lorenz system, where we bound extreme values of various monomials of the coordinates over the global attractor. In the second example we bound extreme values over a chaotic saddle in a nine-mode truncation of fluid dynamics that displays long-lived chaotic transients. The third example has two locally stable limit cycles, each with its own basin of attraction, and we apply our non-global framework to construct bounds for one basin that do not apply to the other. For each example we compute Lyapunov functions of polynomial degrees up to at least eight. In cases where we can judge the sharpness of our bounds, they are sharp to at least three digits when the polynomial degree is at least four or six.
\end{abstract}

\section{Introduction}

In many complicated dynamical systems it is desirable to predict the magnitudes of extreme events---for instance, the greatest instantaneous force applied by a turbulent fluid flow, or the maximum height of a rogue wave. The present work considers extreme events in deterministic systems at late times, as opposed to transient or stochastic behavior. In particular we consider dynamical systems governed by differential equations, especially those with complicated invariant sets such as chaotic attractors or saddles, and we bound the values that quantities of interest can assume on these invariant sets.

When solutions of a differential equation cannot be characterized exactly, a common way to estimate their eventual behavior is to find attracting sets---subsets of state space that attract all bounded sets of initial conditions. The minimal attracting set, which is contained in all others, defines the \emph{global attractor} of the dynamical system~\cite{Robinson2001}. All forward-time limit points are contained in the global attractor. Thus, extreme values of any given quantity at late times are bounded by its extrema over the global attractor, which in turn are bounded by extrema over any other attracting set. Here we construct attracting sets and estimate extrema over them. Various attracting sets can be found using Lyapunov functions, as described below, but generally there exist an infinite number of such functions. A typical approach is to first construct one or several of the simplest possible Lyapunov functions, often quadratic functions, and then use them to estimate properties of the global attractor. In the present work we combine these two steps. For each quantity whose extrema over the global attractor one wants to bound, the construction of a corresponding Lyapunov function that implies the sharpest possible bound is posed as an optimization problem.

Another way to characterize the global attractor is as the maximal compact invariant set~\cite{Robinson2001}, where by definition invariant sets map to themselves under both forward-time and backward-time dynamics. In general a global attractor is different from an attractor, which is typically defined as a minimal \emph{locally} attracting set. An attractor can fail to be a global attractor by not attracting all bounded sets, and a global attracting can fail to be an attractor by not being minimal among locally attracting sets. Since a global attractor is comprised of all compact invariant sets, it will include all attractors but also unstable structures such as fixed points and period orbits, as well as their unstable manifolds.

Our general approach can be carried out analytically or numerically for broad classes of ordinary and partial differential equations (ODEs and PDEs). In practice, however, the construction of optimal Lyapunov functions can be prohibitively difficult, even with computer assistance. An important exception is the case of ODE systems with polynomial right-hand sides. As long as the dimensions and polynomial degrees of such systems are not too large, optimal Lyapunov functions can be constructed computationally using methods of polynomial optimization.

The reason the polynomial ODE case is computationally tractable is that the sufficient conditions we use to construct bounds become, in this case, nonnegativity constraints on polynomial expressions. We enforce this nonnegativity by requiring the expressions to admit decompositions as sums of squares of other polynomials. The optimization over Lyapunov functions subject to these sum-of-squares (SOS) constraints can be carried out numerically after being recast as a semidefinite program (SDP)---a standard type of convex optimization problem \cite{Boyd2004}. The use of SDPs to solve SOS optimization problems was proposed two decades ago \cite{Parrilo2000, Nesterov2000, Lasserre2001} and has found numerous applications in the study and control of ODEs. These applications include the construction of Lyapunov functions to show that a particular solution is attracting and approximate its basin of attraction \cite{Parrilo2000, Papachristodoulou2002, Tan2008, Chesi2011, Henrion2014, Anderson2015a, Valmorbida2017, Drummond2018}, as well as related methods for bounding infinite-time averages \cite{Chernyshenko2014, Fantuzzi2016, Goluskin2018}. Some progress has been made applying SOS methods to nonlinear PDEs also \cite{Goulart2012, Huang2015, Goluskin2019a, Korda2018, Marx2020, Fuentes2019}. These works are similar in spirit to our present method, which was suggested but not applied in \cite{Goluskin2019a}, but they are not the same. As far as we know, the present study is the first to construct Lyapunov functions by optimizing the bounds that they imply.

\Cref{sec: framework} formulates a convex optimization framework for bounding extreme values on global attractors, as well as a modified framework for bounding extrema over sets that attract locally but not globally. Computationally tractable versions of each framework based on SOS constraints are then given for the case of polynomial ODEs. \Cref{sec: ex} reports bounds computed in this way for three polynomial ODEs. For an example with a chaotic attractor, in \cref{sec: Lorenz} we consider the Lorenz system and compute bounds on various monomials of the coordinates. For an example of a chaotic saddle with long transients, in \cref{sec: transient} we study a nine-dimensional ODE that models fluid dynamics in the regime of transition to turbulence \cite{Moehlis2004}. \Cref{sec: multi} concerns an ODE with two locally attracting limit cycles; we test the non-global version of our framework by computing bounds that apply to one limit cycle but not the other. In the first and third example, extreme values can be estimated precisely by searching over particular trajectories, so we can compare our computed bounds to these estimates. In each such case, our computations appear to produce arbitrarily sharp bounds as the polynomial degree of Lyapunov functions is raised. \Cref{sec: con} offers conclusions and open questions.

\section{\label{sec: framework}Constructing bounds using convex optimization}

Consider an autonomous ODE system,
\beq
\ddt\x(t)=\f(\x(t)), \qquad \x(0)=\x_0,
\label{eq: ode}
\eeq
that is well-posed for any initial condition $\x_0\in\R^n$. Assume that all trajectories eventually remain in a bounded subset of $\R^n$, as often can be proved by Lyapunov function methods similar to those used in this work, and that all trajectories are continuously differentiable. The bounding conditions presented below can be generalized beyond ODEs on $\R^n$, including to parabolic PDEs on Banach spaces as in \S2.2 of \cite{Goluskin2019a}. We speak in terms of ODEs here, for simplicity and because this is the case that we can tackle computationally.

The global attractor $\cA$ of \eqref{eq: ode} can be defined as the maximal compact subset of $\R^n$ that is invariant under the dynamics, or equivalently as the minimal set that attracts all initial conditions in every bounded subset of $\R^n$ \cite{Robinson2001}. Let $\Phi:\R^n\to\R$ denote a quantity of interest for solutions to \eqref{eq: ode}. We seek the maximum and minimum values of $\Phi(\x)$ over the global attractor $\cA$,
\begin{align}
\Phi_\cA^+ &:= \max_\cA\Phi(\x), &
\Phi_\cA^- &:= \min_\cA\Phi(\x).
\label{eq: Phi_A}
\end{align}
In particular, our objective is to compute upper bounds on $\Phi_\cA^+$ and lower bounds on $\Phi_\cA^-$. Such results also bound values of $\Phi$ along trajectories $\x(t)$ at sufficiently late times since forward-time limit points lie in the global attractor. That is,
\beq
\Phi_\cA^- ~\le~ \inf_{\x_0\in\R^n}\liminf_{t\to\infty}\Phi(\x(t))
~\le~ \sup_{\x_0\in\R^n}\limsup_{t\to\infty}\Phi(\x(t)) ~\le~ \Phi_\cA^+.
\eeq
It suffices to discuss upper bounds since lower bounds on $\Phi_\cA^-$ are equivalent to upper bounds on $(-\Phi)_\cA^+$. 

In systems with multiple basins of attraction, often one is interested only in the dynamics in a particular basin. Extrema over the global attractor $\cA$ might not be useful in such cases since $\cA$ includes points from every basin. Instead, suppose all trajectories of interest eventually remain in some set $X$. We can define a set $\cA_X$ that is like the global attractor for $X$ alone. That is, $\cA_X$ is the minimal set that attracts all bounded subsets of $X$. In addition to a method for bounding the global maximum $\Phi_\cA^+$, \cref{sec: suff} below gives a method for bounding the non-global maximum
\beq
\Phi_{\cA_X}^+ := \max_{\cA_X}\Phi(\x).
\label{eq: Phi_AB}
\eeq

\subsection{\label{sec: suff}Bounds implied by optimal Lyapunov functions}

We bound $\Phi_\cA^+$ above by constructing attracting sets. Every attracting set contains $\cA$, so the maximum of $\Phi$ over any attracting set is an upper bound on $\Phi_\cA^+$. Like many authors, we find attracting sets using functions $V:\R^n\to\R$, where $V$ is in the set $\mathcal C^1$ of continuously differentiable functions on $\R^n$. We refer to these $V$ as Lyapunov functions, although they need not have properties that sometimes define Lyapunov functions, such as boundedness below or polynomial growth as $|\x|\to\infty$. We aim to show that a sublevel set of $V$,
\beq
\Omega^C_V:=\{\x\in\R^n : V(\x)\le C\},
\eeq
is attracting for some fixed $C$. A sufficient condition for $\Omega^C_V$ to be attracting is that
\beq
\lambda\,\f(\x)\cdot\nabla V(\x) \le C - V(\x)
\label{eq: lyap}
\eeq
throughout $\R^n$ for some $\lambda>0$. To see why the above suffices, note that $\ddt V(\x(t)) = \f(\x(t))\cdot\nabla V(x(t))$ along all trajectories of \eqref{eq: ode}, so \cref{eq: lyap} implies
\beq
\ddt\left[V(\x(t))-C\right] \le -\frac1\lambda\left[V(\x(t))-C\right].
\label{eq: diff ineq}
\eeq
Integrating the differential inequality~\cref{eq: diff ineq} using an integrating factor gives $V(\x(t))-C\le e^{-t/\lambda}[V(\x_0)-C]$. This upper bound decreases to zero as $t\to\infty$, and $V$ and trajectories are continuous, so the set $\Omega^C_V$ attracts all bounded trajectories.

The attracting set $\Omega_V^C$ contains $\cA$ by definition, so
\beq
\Phi_\cA^+ \le \max_{\x\in\Omega_V^C} \Phi(\x).
\label{eq: Omega max}
\eeq
The right-hand maximum in \cref{eq: Omega max} can be prohibitively difficult to evaluate for complicated $V$. We avoid this difficulty by adding a second constraint on $V$, requiring not only \eqref{eq: lyap} but also $\Phi(\x)\le V(\x)$ throughout $\R^n$. Then the right-hand maximum in \eqref{eq: Omega max} is bounded above by $\max_{\x\in\Omega_V^C} V(\x)$, whose value is simply $C$, so we arrive at the upper bound $\Phi_\cA^+\le C$. Optimizing this bound over all choices of $\lambda$ and $V$ gives
\beq
\Phi_\cA^+ \le \inf_{\substack{\lambda>\R\\V\in\,\mathcal C^1}}C \quad s.t. \quad 
\begin{array}[t]{rl}
V(\x)-\Phi(\x) & \ge0 \quad\forall \x\in\R^n, \\ 
	C - V(\x) - \lambda\,\f(\x)\cdot\nabla V(\x) & \ge 0 \quad\forall \x\in\R^n.
\end{array}
\label{eq: opt general}
\eeq
For various computational examples in \cref{sec: ex}, the left-hand inequality in \eqref{eq: opt general} appears to be an equality, meaning that our bounding conditions can give arbitrarily sharp bounds on $\Phi_\cA^+$. It remains an open challenge to prove this equality in general under suitable conditions on $\f$ and $\Phi$.

The optimization \eqref{eq: opt general} is convex in $V$ for each fixed value of $\lambda$, and the same is true when the convex space $\mathcal C^1$ is replaced by any convex subspace. This convexity makes it tractable to optimize $V$ computationally in certain cases, at least over finite-dimensional subsets of $\mathcal C^1$, which is why we favor the inequalities in \eqref{eq: opt general} over other sufficient conditions for $\Phi_\cA^+\le C$. Many previous authors have taken a different approach, choosing a particular function or simple ansatz for $V$ at the start of their analyses. With $V$ so fixed, one can use more complicated sufficient conditions that might give attracting sets $\Omega_V^C$ with smaller values of $C$ than can be obtained using \eqref{eq: lyap}.\footnote{One weaker sufficient condition is to let $C$ be the maximum of $V$ on the set where $\f\cdot\nabla V$ vanishes \cite{Krishchenko2005}. This amounts to imposing \eqref{eq: lyap} on that set instead of on $\R^n$; the value $C$ is attained at a stationary point of the Lagrangian $V+\lambda\,\f\cdot\nabla V$, where here $\lambda$ is a Lagrange multiplier. Additionally, constraints on $V$ can be restricted to subsets of $\R^n$ already known to be attracting, and various attracting sets can be intersected to produce a smaller attracting set.}
Nonetheless, there are many examples in the literature where the best bounds on $\Phi_\cA^+$ are not close to being sharp. We propose that better bounds can be obtained by considering larger classes of Lyapunov functions, even with our less powerful sufficient conditions. This is borne out by the computational examples of \cref{sec: ex}.

Any $V$ that gives a bound $\Phi_\cA^+\le C$ also gives information about where on the global attractor the value of $\Phi$ can be close to $C$. Let $\x^\eps$ denote any point on $\cA$ where \mbox{$\Phi(\x^\eps)+\eps\ge C$}. (If the bound $C$ is not sharp, such a point exists only when $\eps$ is sufficiently large.) The constraints on $V$ in \cref{eq: opt general} and the fact that $\x^\eps$ is in $\Omega_V^C$ imply that $C-\eps\le V\le C$ and $\f\cdot\nabla V\le\eps/\lambda$ at $\x^\eps$. For fixed $(V,C,\lambda,\eps$), these inequalities define a subset of $\R^n$ in which any near-optimal $\x^\eps$ must lie. Such subsets could be used to identify parts of state space in which extreme behavior occurs, similar to what is done for other types of extreme behavior in \cite{Tobasco2018} and \cite{Fantuzzi2019}, but we do not pursue this idea here.

A small modification to the framework \eqref{eq: opt general} can be used to bound maxima over $\cA_X$, the minimal set that attracts all bounded subsets of $X$. The same arguments leading to \eqref{eq: opt general} give an optimization problem for non-global bounds:
\beq
\Phi_{\cA_X}^+ \le \inf_{\substack{\lambda>\R\\V\in\,\mathcal C^1(X)}}C \quad s.t. \quad 
\begin{array}[t]{rl}
V(\x)-\Phi(\x) & \ge0 \quad\forall \x\in X, \\ 
	C - V(\x) - \lambda\,\f(\x)\cdot\nabla V(\x) & \ge 0 \quad\forall \x\in X.
\end{array}
\label{eq: opt X}
\eeq
It some cases it may be difficult to find a suitable choice of $X$. For instance, one may be interested in a particular basin of attraction, but typically such a basin is not known exactly and may be a fractal. Then the challenge is to find a choice of $X$ that is small enough to lie inside the basin and large enough to attract every trajectory inside the basin.

\subsection{\label{sec: sos}Computing bounds by polynomial optimization}

Optimization over Lyapunov functions as in \eqref{eq: opt general} can be carried out by methods of polynomial optimization if the ODE right-hand side $\f(\x)$ and quantity of interest $\Phi(\x)$ are both polynomials. Henceforth we assume this is the case, and we optimize Lyapunov functions not over all of $\mathcal C^1$ but over the the set of real polynomials in $n$ variables up to a specified degree $d$---that is, the set $\R[\x]_{n,d}$. The inequalities in \eqref{eq: opt general} then require nonnegativity of two multivariable polynomials on $\R^n$. Deciding whether a polynomial is nonnegative has NP-hard computational complexity unless $n$ or $d$ is small \cite{Murty1987}, and we want to optimize among higher-degree Lyapunov functions for which such computations would be intractable. Thus we employ a standard SOS relaxation, replacing nonnegativity of a polynomial with the generally stronger constraint that the polynomial can be represented as a sum of squares of other polynomials \cite{Parrilo2013a}. The resulting SOS optimization problem is
\beq
\Phi_\cA^+ \le \inf_{\lambda>\R}~\inf_{V\in\R[\x]_{n,d}}C \quad s.t. \quad 
\begin{array}[t]{rl}
V(\x)-\Phi(\x) & \in\Sigma_n,
	 \\ C - V(\x) - \lambda\,\f(\x)\cdot\nabla V(\x) & \in\Sigma_n,
\end{array}
\label{eq: opt sos}
\eeq
where $\Sigma_n$ denotes the set of SOS polynomials in $n$ variables. The optimization in \eqref{eq: opt sos}, relative to that in \eqref{eq: opt general}, has a smaller set of admissible $V$ and so gives an upper bound that is at least as large. This bound improves or remains unchanged as the polynomial degree $d$ of $V$ is raised. Convergence as $d\to\infty$ is discussed below.

The inner minimization in \eqref{eq: opt sos} can be reformulated as a semidefinite program (SDP) because the SOS constraints and optimization objective are linear in the tunable variables, which are $C$ and the coefficients of $V$ \cite{Parrilo2013a}. This linearity is why we do not optimize over $\lambda$ simultaneously, instead tuning $C$ and $V$ with fixed $\lambda$. The outer minimization might be difficult in some cases since the dependence of the inner minimum on $\lambda$ need not be convex or even continuous. It is only a one-dimensional search, however, and at least for the examples of \cref{sec: ex} we find simple dependence on $\lambda$.

The non-global bounding formulation \eqref{eq: opt X} also can be altered so that nonnegativity on $X$ is enforced by SOS constraints that imply nonnegativity on $X$ but not on $\R^n$. We assume the set $X$ is semialgebraic, meaning it can be specified by a finite number of polynomial inequalities and equalities. For convenience, suppose it can be specified by a single inequality, so there exists a polynomial $g$ such that $g(\x)\ge0$ if and only if $\x\in X$. A sufficient condition for any polynomial $p$ to be nonnegative on $X$, without necessarily being nonnegative outside of $X$, is that there exists an SOS polynomial $s$ such that $p-gs$ also is SOS. This standard approach, which often is called the $\mathcal S$-procedure \cite{Tan2008} or a weighted SOS condition \cite{Lasserre2015}, readily gives an SOS relaxation of \eqref{eq: opt X}:
\beq
\Phi_{\cA_X}^+ \le \inf_{\lambda>\R}~\inf_{\substack{V\in\R[\x]_{n,d}\\ \hspace{3pt}s_1\in\R[\x]_{n,d_1}\\\hspace{3pt}s_2\in\R[\x]_{n,d_2}}}C \quad s.t. \quad 
\begin{array}[t]{rl}
V(\x)-\Phi(\x)-g(\x)s_1(\x) & \in\Sigma_n,
	 \\ C - V(\x) - \lambda\,\f(\x)\cdot\nabla V(\x) - g(\x)s_2(\x) & \in\Sigma_n,
	 \\ s_1(\x),~s_2(\x) & \in\Sigma_n.
\end{array}
\label{eq: opt X sos}
\eeq
As in \eqref{eq: opt sos}, the inner minimization problem can be reformulated as an SDP and solved computationally.

For a general semialgebraic set $X$ that is specified by inequalities $g_j(\x)\ge0$ and equalities $h_k(\x)=0$, each $g(\x)s_i(\x)$ term in \eqref{eq: opt X sos} can be replaced by a term of the form $\sum_{j=1}^Jg_j(\x)s_j(\x)+\sum_{k=1}^Kh_k(\x)r_k(\x)$, where the $s_j$ are SOS but the $r_k$ are arbitrary polynomials. If the specification of $X$ includes an inequality of the form $g_j(\x)=R^2-|\x|^2\ge0$ for some constant $R$ (which for compact $X$ always can be added without changing the set), then Putinar's Positivstellensatz theorem \cite{Putinar1993} implies that the optimum of the SOS problem converges to that of the original optimization \eqref{eq: opt X} as the polynomial degrees approach infinity \cite{Lasserre2015}. For the global formulation, we cannot say in general that the optimum of \eqref{eq: opt sos} converges to that of \eqref{eq: opt general} as $d\to\infty$. If such a convergence guarantee is important, one option is to find $R$ such that $g(\x)=R^2-|\x|^2\ge0$ on the global attractor $\cA$, and then solve the right-hand optimization problem in \eqref{eq: opt X sos} instead of \eqref{eq: opt sos}. In practice, no such modification may be required; in the example of \cref{sec: Lorenz}, the global SOS formulation \eqref{eq: opt sos} provides very sharp bounds on $\Phi_\cA^+$ as $d$ is raised.

\subsection{\label{sec: structure}Exploiting structure}

To reduce computational cost and improve numerical conditioning when solving \eqref{eq: opt sos} or \eqref{eq: opt X sos}, often one can impose some structure \emph{a priori} on the polynomial ans\"atze of $V$ without changing the optimal bounds. In the examples of \cref{sec: ex}, we impose structure in two ways.

The first source of structure is symmetry. Let $\mathcal T:\R^n\to\R^n$ be an invertible linear transformation that generates a finite group, meaning that composing $\mathcal T$ a finite number of times gives the identity. Suppose the ODE and the quantity to be bounded both are invariant under $\mathcal T$, meaning $\ddt\mathcal T\x=\f(\mathcal T\x)$ and $\Phi(\mathcal T\x)=\Phi(\x)$. (Each ODE in \cref{sec: ex} has such a symmetry, which is shared by some of the $\Phi$ we bound but not all.) Imposing the invariance $V(\mathcal T\x)=V(\x)$ does not change the optimum of \eqref{eq: opt sos}. We omit the proof of this statement because it is closely analogous to the proof of Proposition A.1 in \cite{Goluskin2019a} or Theorem 2 in \cite{Lakshmi2020}. Essentially, any non-symmetric $V$ can be averaged over the group orbit of $\mathcal T$ to obtain a symmetrized $V$ that gives the same bounds. In the non-global formulation \eqref{eq: opt X sos}, if the function $g$ defining $X$ also is invariant under $\mathcal T$, then we can impose $\mathcal T$-invariance on $V$ and each $s_i$ without changing the optimal bound. With $\mathcal T$-invariance imposed on $V$ in \cref{eq: opt sos}, or on $V$ and the $s_i$ in \cref{eq: opt X sos}, the expressions that are constrained to be SOS also are $\mathcal T$-invariant. This can be exploited when reformulating the SOS optimizations as SDPs \cite{Parrilo2013a}, leading to smaller and better-conditioned SDPs. 

There is a second source of structure when the \emph{global} bounding formulation~\eqref{eq: opt sos} is applied to an ODE whose right-hand side $\f$ has highest-degree terms of even polynomial degree. In the typical situation where $\deg(V)>\deg(\Phi)$, the first and second constraints in~\eqref{eq: opt sos} require that the highest-degree terms in $V$ and $\f\cdot\nabla V$, respectively, both are of even degree. (There is no such requirement in the non-global formulation~\eqref{eq: opt X sos} since the highest-degree terms can instead come from $gs_1$ and $gs_2$.) When $\deg(f)$ and $\deg(V)$ both are even, one generally expects $\deg(\f\cdot\nabla V)$ to be odd. However, the degree of $\f\cdot\nabla V$ can be reduced by one if the highest-degree terms in $V$ cancel in the expression $\f\cdot\nabla V$. Requiring this cancellation gives constraints on the coefficients of $V$ that can be deduced \emph{a priori}. We do so when applying the global formulation in the examples of \cref{sec: Lorenz} and \cref{sec: transient}, where $\f$ is quadratic.

\section{\label{sec: ex}Computational examples}

Each of the following three subsections reports bounds that we have computed for an ODE example. In the first two examples we compute bounds over global attractors by solving \eqref{eq: opt sos}, and in the last example we compute non-global bounds by solving \eqref{eq: opt X sos}. In all computations, the parser YALMIP \cite{Lofberg2004, Lofberg2009} (version R20190425) was used to translate SOS optimizations into SDPs, which were solved using MOSEK 9.0. In all cases MOSEK converged with relative infeasibilities below $5\cdot10^{-7}$.

The numerical computations are subject to rounding error, so the reported bounds are not rigorous to the standards of a computer-assisted proof. Rather, much as numerical integration approximates a particular solution with uncontrolled but often small rounding error, our computations approximate global statements about dynamics. It is possible to make SOS computations rigorous using interval arithmetic, as was done in \cite{Goluskin2018} for the SDP solver but not the parser, however we do not do so here.

\subsection{\label{sec: Lorenz}The Lorenz system}

To test the quality of bounds computed using \eqref{eq: opt sos} for a system with a chaotic attractor, we consider the Lorenz equations \cite{Lorenz1963}, for which the components of the generic ODE \eqref{eq: ode} are
\beq
\x = (x,y,z), \quad 
\f = (-\sigma x + \sigma y \;,\; rx - y - xz \;,\; - \beta z + xy).
\label{eq: Lorenz}
\eeq
We fix the parameters to their standard chaotic values $(\beta,\sigma,r)=(8/3,10,28)$, at which there exists a strange attractor to which almost every trajectory tends \cite{Tucker1999}. Invariant structures embedded in the strange attractor include an equilibrium at the origin and an infinite number of periodic orbits, as well as their unstable manifolds \cite{Sparrow1982, Doedel2015}. The global attractor $\cA$ includes all such structures, as well as the invariant structures which are not part of the strange attractor: the two equilibria $\x_\pm = (\pm6\sqrt{2},\pm6\sqrt{2},27)$ and their unstable manifolds. As described below, the $\Phi$ we bound turn out to be extremized on the strange attractor, rather than on the unstable manifolds of $\x_\pm$, so the sharp bounds we compute on $\Phi_\cA^+$ and $\Phi_\cA^-$ are also sharp bounds for the strange attractor alone.

The quantities we have bounded are various monomials $\Phi=x^ly^mz^n$ up to cubic degree, although our framework applies just as easily to more general polynomials. Bounds were computed by solving \eqref{eq: opt sos} as described above---sweeping through $\lambda$ and solving the inner minimization as an SDP---with $V(x,y,z)$ of polynomial degrees up to 8. Numerous authors have found that, for SDPs coming from dynamically motivated SOS optimizations, numerical conditioning is improved if ODE variables are rescaled so that the relevant dynamics lie approximately in the cube $[-1,1]^n$. Here such rescaling is again crucial. We solved \eqref{eq: opt sos} with the Lorenz equations rescaled by $(x,y,z)\mapsto25(x,y,z)$ as in \cite{Goluskin2018}, then the computed bounds were converted back to the original scaling.

Rather than using fully general polynomial ans\"atze for $V$, we can impose some structure without changing the optimal bounds, as described in \cref{sec: structure}. First, since $\f$ is quadratic and $\deg(V)$ and $\deg(\f\cdot\nabla V)$ both must be even, the highest-degree terms in $V$ must cancel in $\f\cdot\nabla V$. Such cancellation occurs if and only if all highest-degree terms of $V$ take the form $x^p(y^2+z^2)^q$ \cite{Swinnerton-Dyer2001, Goluskin2018}, so we impose this on $V$. Second, for $\Phi$ that are invariant under the symmetry $(x,y)\mapsto(-x,-y)$ of the Lorenz equations, we can impose the same invariance on $V$ without changing the optimal bounds.

In order to judge the sharpness of our bounds we have sought extreme values of each $\Phi$ among particular trajectories of the Lorenz equations. Such a search might be impossible in more complicated systems, which is one motivation for our bounding approach, but it is possible here. Trajectories we examined include numerical integrations beginning from random initial conditions (with initial transients removed), numerical integration approximating the one-dimensional unstable manifold of the origin, and the many periodic orbits computed by Viswanth \cite{Viswanath2003, Viswanath2004}. It is on the origin or its unstable manifold that we find the extreme values of each $\Phi$, and the closeness of these values to our computed bounds suggests that they are indeed global extrema. Numerical integration from random initial conditions does not give very good approximations to these various extrema; integrating for $10^7$ fourth-order Runge--Kutta time steps of size 0.005 gives $\Phi$ values that share one or two digits with the true extrema but no more.

\begin{figure}[tp]
\begin{center}
\includegraphics[width=220pt,trim={20pt 2pt 46pt 10pt},clip]{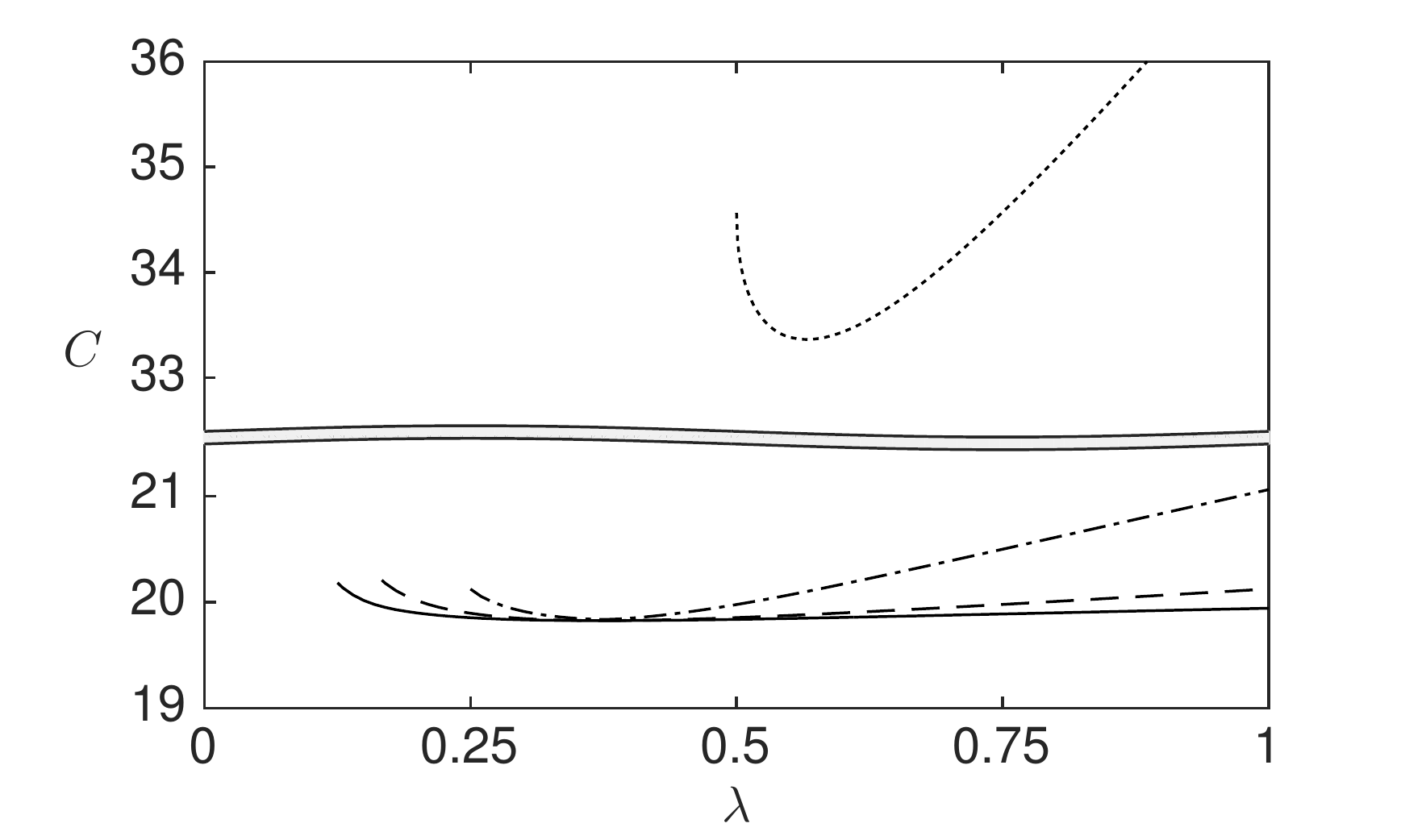}
\end{center}
\caption{\label{fig: lambda}Upper bounds $C$ on the maximum of $x$ over the global attractor of the Lorenz equations at the standard parameters. The bounds are optima of the inner minimization in \eqref{eq: opt sos} for various $\lambda$ and $V$ of polynomial degrees 2 (\dottedrule), 4 (\dashdottedrule), 6 (\longdashedrule), and 8 (\solidrule).}
\end{figure}

Bounds produced by the inner minimization problem in \eqref{eq: opt sos} depend on $\lambda$ and the degree of $V$ in similar ways for all $\Phi$ that we considered for the Lorenz equations. As a typical instance, \cref{fig: lambda} shows this dependence for upper bounds on the maximum of $x$ over $\cA$. For all $\Phi$, our computations give finite bounds with degree-$d$ Lyapunov functions when $\lambda\in[1/d,\infty)$, and the bounds are convex in $\lambda$ on these intervals. Thus it is not hard to optimize $\lambda$ over these intervals, which apparently suffices to give arbitrarily sharp bounds as the degree of $V$ is raised. However, this dependence of $C$ on $\lambda$ is particular to the Lorenz equations.\footnote{To see why finite bounds for the Lorenz equations require $\lambda\ge1/d$, note that terms in $V$ of the form $c_1x^d$, $c_2y^d$, and $c_3z^d$ produce terms in $-(V+\lambda\,\f\cdot \nabla V)$ of the form $c_1(-1+\sigma\lambda d)x^d$, $c_2(-1+\lambda d)y^d$, and $c_3(-1+\beta\lambda d)z^d$, respectively. The latter three coefficients must be nonnegative in order for the second constraint in \eqref{eq: opt sos} to hold. At the standard parameters this requires $\lambda\ge1/d$ if $c_1,c_2,c_3>0$. In almost all cases we obtain finite bounds only when $c_1,c_2,c_3>0$, and thus only when $\lambda\ge 1/d$ also. An exception is the lower bound on $z$ with quadratic $V$, where finite bounds are possible with $c_2=0$ and $\lambda=3/8\notin[1/2,\infty)$.}

\begin{figure}[tp]
\begin{center}
\includegraphics[width=210pt,trim={12pt 2pt 36pt 20pt},clip]{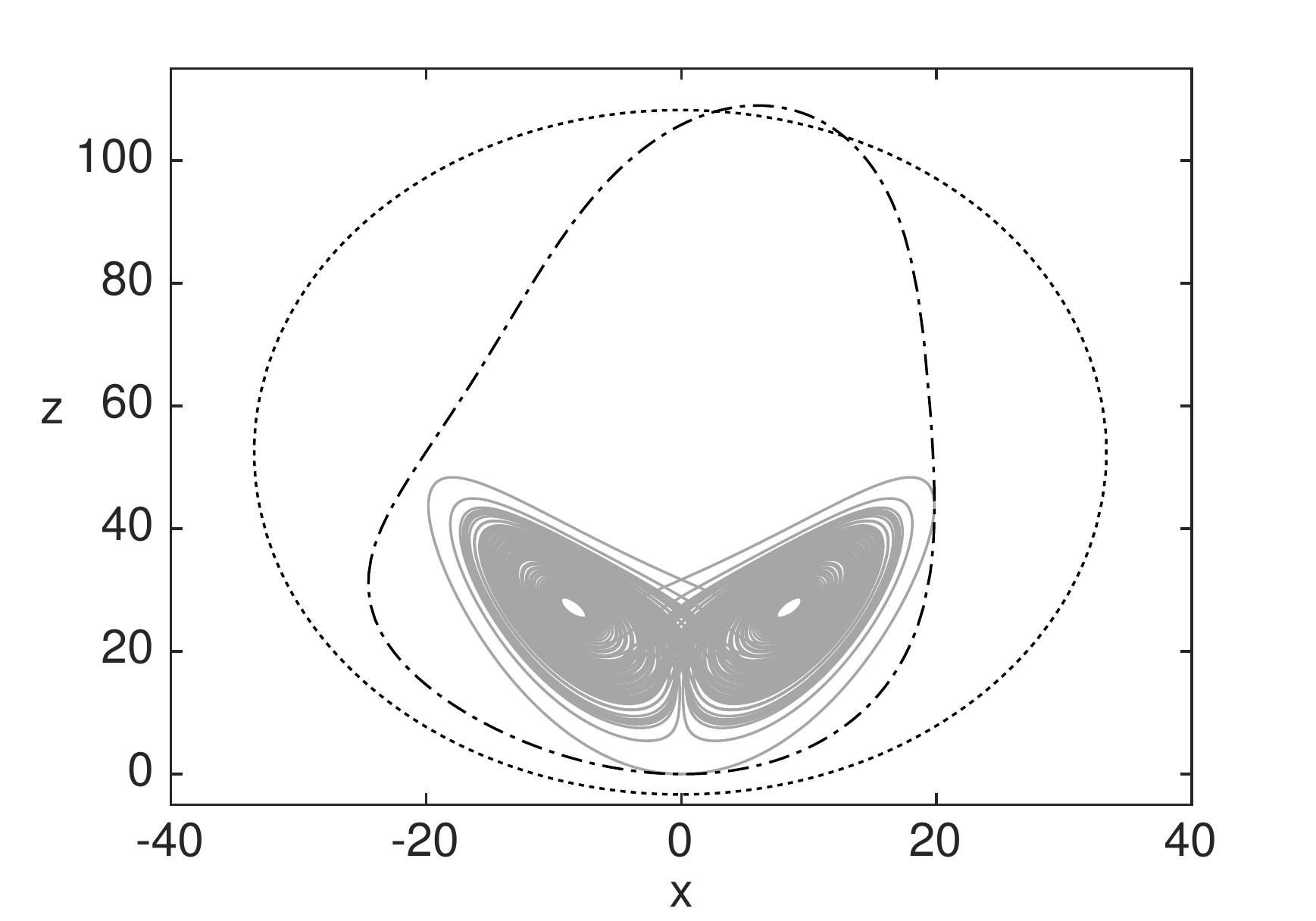}
\end{center}
\caption{\label{fig: boundaries}Attracting sets $\Omega_V^C$ that give our best upper bounds $C$ on the maximum of $x$ over the global attractor of the Lorenz equations, for Lyapunov functions $V$ of degree 2 (\dottedrule) and 4 (\dashdottedrule). Each pair $(V,C)$ solves the inner minimization in \eqref{eq: opt sos} with $\lambda=0.5659$ and $\lambda=0.3743$, respectively. The plotted curves are boundaries of the projections of $\Omega_V^C$ onto the $xz$-plane. Also shown are numerically integrated trajectories (${\color{mygray}\solidrule}$) starting along each half of the origin's unstable manifold.}
\end{figure}

As an example of the attracting sets $\Omega_V^C$ that give the upper bounds on $x$ reported in \cref{fig: lambda}, let us consider quadratic and quartic $V$ at the approximately optimal values of $\lambda=0.5659$ and $\lambda=0.3743$, respectively. \Cref{fig: boundaries} shows the quadratic and quartic attracting sets $\Omega_V^C$, where $V$ and $C$ solve the inner minimization in \eqref{eq: opt sos} at the specified $\lambda$ values. Also shown in \cref{fig: boundaries} is a numerical approximation to the unstable manifold of the origin, which is part of the strange attractor.

\begin{figure}[tp]
\begin{center}
\includegraphics[width=240pt,trim={14pt 2pt 32pt 20pt},clip]{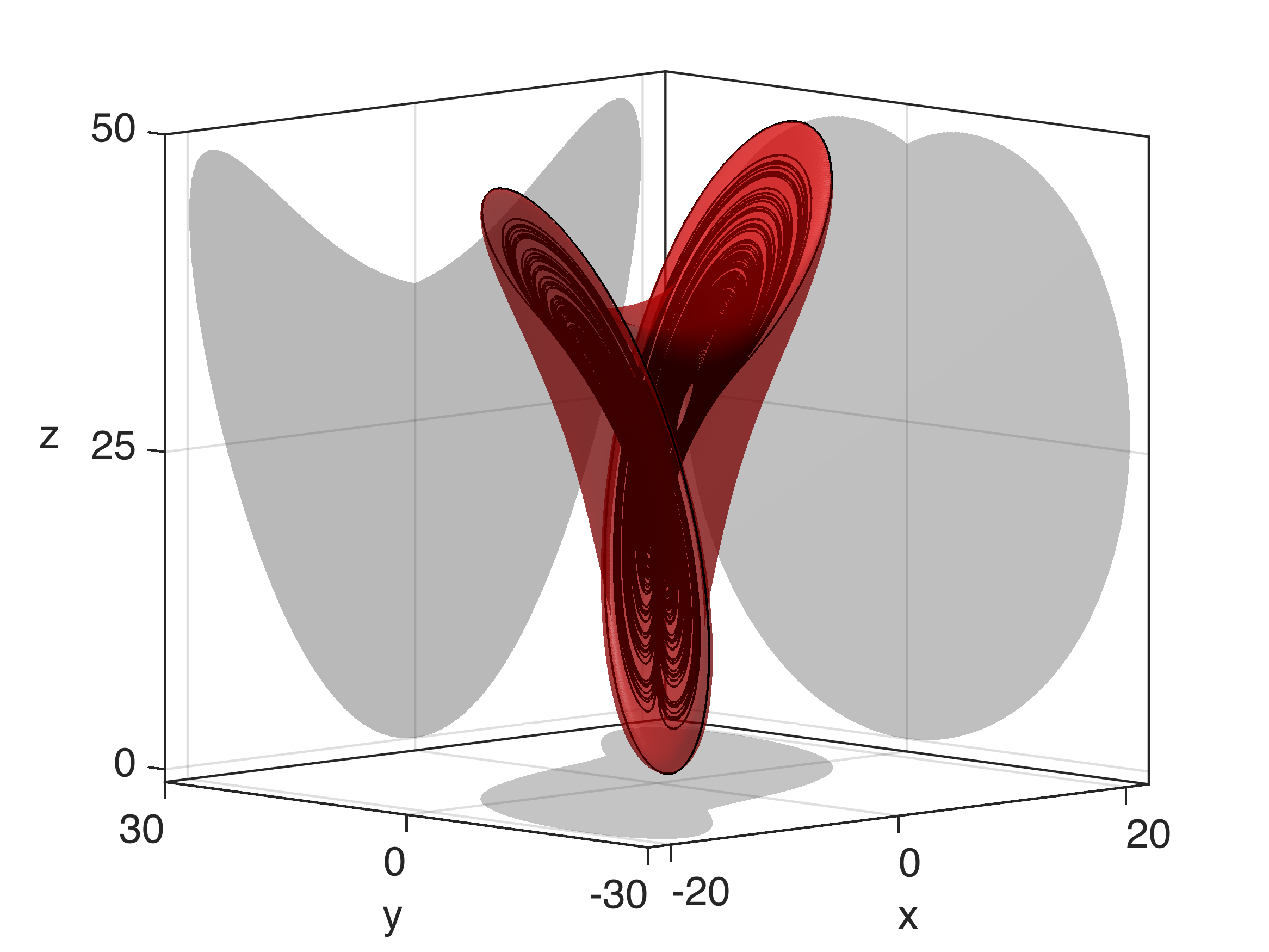}
\end{center}
\caption{\label{fig: 3D}Attracting set $\Omega_V^C$ for $\lambda=3/8$ and the degree-8 $V$ that optimizes the resulting bound $C$ on the maximum of $z$ over the global attractor of the Lorenz equations. Also shown are numerically integrated trajectories starting along each half of the origin's unstable manifold.}
\end{figure}

Attracting sets constructed by solving \eqref{eq: opt sos} do not necessary localize the global attractor well since this is not the optimization objective. The quartic attracting set in \cref{fig: boundaries} provides a very good bound on $x$ but a poor localization of the global attractor as a whole. Nonetheless, some other attracting sets constructed by solving \eqref{eq: opt sos} do localize the global attractor well. \Cref{fig: 3D} shows one such $\Omega_V^C$, which was constructed by optimizing upper bounds on $z$ using $V$ of degree 8. This may be the smallest attracting set that has been reported for the Lorenz equations. It would be harder to minimize the volume of $\Omega_V^C$ directly because this volume does not have convex dependence on the coefficients of $V$.

For all monomials $x^ly^mz^n$ considered here, we have computed bounds on extreme values over the global attractor of the Lorenz equations that are either exactly sharp or very close to being so. For convenience we normalize our bounds by each monomial's value at the nonzero fixed points $\x_\pm$, meaning that we report results for
\beq
\Phi = \frac{x^ly^mz^n}{|x^ly^mz^n|_{\x_\pm}} = \frac{x^ly^mz^n}{(6\sqrt2)^{l+m}27^n},
\label{eq: phi norm}
\eeq
although the bounds were computed using unnormalized monomials. \Cref{tab: UB} reports upper bounds on these $\Phi$ for all monomials up to cubic degree, along with the apparently maximal values of $\Phi$ found on the unstable manifold of the origin. Quadratic Lyapunov functions, to which many past studies have been confined, do not produce particularly good bounds. On the other hand, $V$ of degree 4 and 6 produce upper bounds for all 13 monomials that are sharp to at least 2 and 4 digits, respectively.

\begin{table}[p]
\caption{\label{tab: UB}Upper bounds on maxima over the global attractor of the Lorenz equations for all normalized monomials \eqref{eq: phi norm} up to cubic degree, computed by solving \eqref{eq: opt sos} with $V$ of degrees 2, 4, and 6. Powers of $x,y,z$ are omitted since their extrema are determined by the extrema of $x,y,z$. Also shown is an approximation of each $\Phi_\cA^+$ from below, found by numerical integration of the unstable manifold of the origin. All values have been rounded to the precision shown. The values of $\lambda$ giving the tabulated bounds are reported by \cref{tab: UB lambda} in \Cref{app: results}.}
\begin{center}
\begin{tabular}{ccccc}
Monomial & \multicolumn{3}{c}{Upper bounds on $\Phi_\cA^+$} 
	& Approximate $\Phi_\cA^+$ \\[2pt] \hline
& $\deg=2$ & $\deg=4$ & $\deg=6$ \\ \cline{2-4}
$x$ & 3.9268 & 2.3378 & 2.3365 & 2.3365 \\ 
$y$ & 3.4081 & 3.2630 & 3.2630 & 3.2630 \\ 
$z$ & 2.1081 & 1.7943 & 1.7912 & 1.7912 \\ 
$xy$ & 10.1143 & 6.8780 & 6.8699 & 6.8698 \\ 
$xz$ & 7.8971 & 3.9948 & 3.9872 & 3.9872 \\ 
$yz$ & 4.6415 & 4.0834 & 4.0832 & 4.0832 \\ 
$x^2y$ && 15.2384 & 15.2288 & 15.2288 \\ 
$x^2z$ && 9.1943 & 9.1619 & 9.1617 \\ 
$xy^2$ && 21.9543 & 21.9483 & 21.9483 \\ 
$xyz$ && 9.4056 & 9.3945 & 9.3944 \\ 
$xz^2$ && 7.0518 & 7.0276 & 7.0276 \\ 
$y^2z$ && 12.2374 & 12.2258 & 12.2258 \\ 
$yz^2$ && 6.1676 & 6.1668 & 6.1668 
\end{tabular}
\end{center}
\end{table}

\begin{table}[p]
\caption{\label{tab: LB}Lower bounds on minima over the global attractor of the Lorenz equations for all normalized symmetric monomials \eqref{eq: phi norm} up to cubic degree, computed by solving \eqref{eq: opt sos} with $V$ of degrees 2, 4, 6, and 8. Also shown is an approximation of each $\Phi_\cA^-$ from above, found by numerical integration of the unstable manifold of the origin. All values have been rounded to the precision shown. The values of $\lambda$ giving the tabulated bounds are reported by \cref{tab: LB lambda} in \Cref{app: results}.}
\begin{center}
\begin{tabular}{cccccc}
Monomial & \multicolumn{4}{c}{Lower bounds on $\Phi_\cA^-$} 
	& Approximate $\Phi_\cA^-$ \\[2pt] \hline
& $\deg=2$ & $\deg=4$ & $\deg=6$ & $\deg=8$ \\ \cline{2-5}
$z$ & 0 &&&& 0 \\
$xy$ & $-10.1143$ & $-1.5644$ & $-0.9048$ & $-0.9043$ & $-0.9042$ \\  
$x^2z$ && $-0.3484$ & $-0.0177$ & $-0.0013$ & 0 \\
$xyz$ && $-2.5369$ & $-1.3920$ & $-1.3914$ & $-1.3914$  \\ 
$y^2z$ && $-0.2898$ & $-0.0309$ & $-0.0061$ & 0 \\
\end{tabular}
\end{center}
\end{table}

As for minima over the global attractor of the Lorenz system, many lower bounds can be anticipated without additional computation. Manifestly nonnegative $\Phi$ such as $x^2$ attain their minima on the equilibrium at the origin. For $\Phi$ that are antisymmetric under the symmetry $(x,y)\mapsto(-x,-y)$ of the Lorenz equations, the upper bound $\Phi_\cA^+\le C$ implies the lower bound $-C\le \Phi_\cA^-$. Thus we compute lower bounds only for symmetric monomials that are not obviously nonnegative. \Cref{tab: LB} reports our lower bounds on such monomials up to cubic degree, computed by solving \eqref{eq: opt sos} for upper bounds on $-x^ly^mz^n$. In all cases, the bounds appear to become sharp as the degree of $V$ is raised.

Complications arise with the lower bounds on $z$, $x^2z$, and $y^2z$ that did not arise with our upper bounds. These three quantities are minimized at the origin, whereas all other extrema we have bounded appear to occur elsewhere on the origin's unstable manifold. With quadratic $V$, the sharp lower bound $0\le z$ can be proved using the sufficient condition~\eqref{eq: lyap} only if $\lambda=3/8$, but a naive search over $\lambda$ may not find this result since other $\lambda$ values smaller than $1/2$ do not give finite bounds.\footnote{With $\lambda=1/\beta$, $C=0$, and $\Phi=-z$, the quadratic Lyapunov function $V=-z+\tfrac{1}{2\sigma}x^2$ satisfies \eqref{eq: lyap} and thus proves the known result $z\ge\tfrac{1}{2\sigma}x^2$. In the past this has been proved by showing that a condition like \eqref{eq: lyap} holds on a compact set already known to be attracting \cite{Sparrow1982, Giacomini1997}, but choosing $\lambda=3/8$ lets \eqref{eq: lyap} hold on all of state space.}
Raising the degree of $V$ to $d\ge4$ removes this difficulty since then the value $\lambda=3/8$ falls in the interval $[1/d,\infty)$ over which the lower bound on $z$ is convex in $\lambda$. The nonnegativity of $z$ on the global attractor implies that $x^2z$ and $y^2z$ are nonnegative also. However, as reflected in \cref{tab: LB}, we have not been able to prove exact lower bounds on $x^2z$ and $y^2z$ directly with the framework of \eqref{eq: opt sos}.

Various bounds on coordinates of the Lorenz equations have appeared in the literature, and bounds on other functions of $(x,y,z)$ can be inferred from previously known attracting sets. The bounds reported in \cref{tab: UB,tab: LB} are sharper than the best results in the literature, except for the already sharp lower bound $0\le z$. For the example of upper bounds on $y$ and $z$, the best prior results we know of are identical to the bounds we report in \cref{tab: UB} for quadratic $V$; both bounds follow from the fact that the cylinder $y^2+(z-r)^2\le\tfrac{\beta^2r^2}{4(\beta-1)}$ is attracting when $\beta\ge2$ \cite{Leonov1987, Doering1995b}. These bounds exceed the true maxima of $y$ and $z$ by more than 4\% and 17\%, respectively, whereas the bounds we compute with quartic $V$ are much sharper. While most authors have considered only quadratic Lyapunov functions for the Lorenz equations, a few have suggested particular quartic functions \cite{Pogromsky2003, Krishchenko2005, Suzuki2008}. None of these quartic functions do as well as our optimized quartic $V$, although the quartic attracting set of \cite{Suzuki2008} implies bounds on $y$ and $z$ that are slightly better than our quadratic-$V$ results.

Some bounds in the literature on the Lorenz equations use analyses more complicated than the sufficient condition \eqref{eq: lyap}. The best prior upper bound on $x$ seems to be that of \cite{Krishchenko2006}, whose approach \cite{Krishchenko1997} is to first use a quadratic Lyapunov function to show that a certain ellipsoid is attracting, and then use $V=|x|$ as a Lyapunov function on that ellipsoid. The resulting bound, normalized according to \eqref{eq: phi norm}, is about 3.180. This is sharper than our quadratic-$V$ bound of 3.9317 but not our quartic-$V$ bound of 2.3378. Similarly, a large number of quadratic Lyapunov functions are constructed in \cite{Robenack2018} using computer algebra, and the implied bounds are stronger than those of a single quadratic $V$ but weaker than those of quartic $V$. These results reflect the fact that the sufficient condition \eqref{eq: lyap} is not the strongest possible. However, they also suggest that inferring the best possible bound from a particular $V$ is not as important as having a computationally tractable way to optimize over $V$ beyond the quadratic case.

\begin{figure}[tp]
\begin{center}
\includegraphics[width=220pt,trim={14pt 0pt 34pt 26pt},clip]{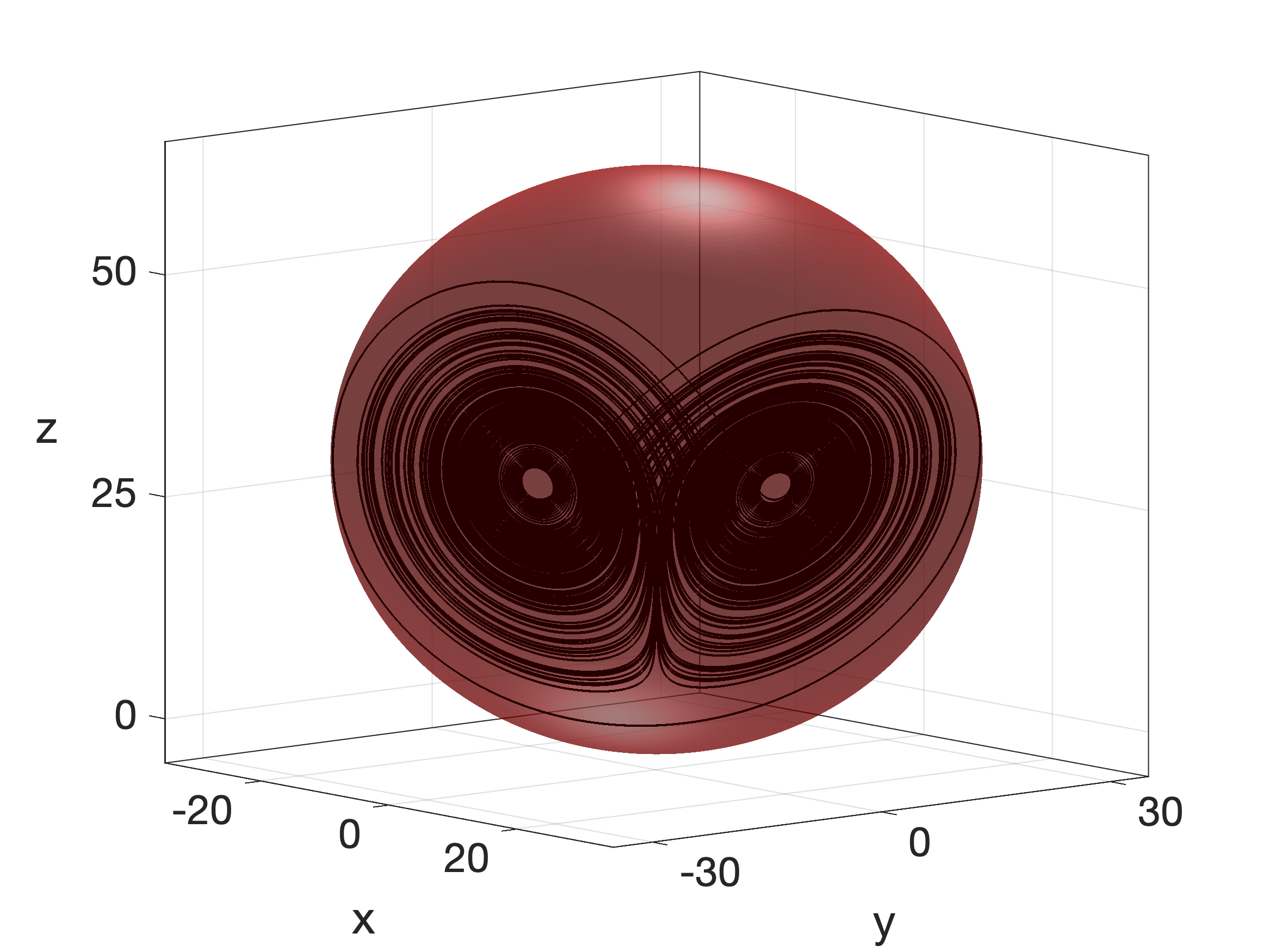}
\end{center}
\caption{\label{fig: sphere}A sphere of nearly minimal radius that contains the global attractor of the Lorenz equations, found using degree-6 $V$ as described in the text. Also shown are numerically integrated trajectories starting along each half of the origin's unstable manifold.}
\end{figure}

As a final example for the Lorenz system, suppose we want to find the smallest possible ball that contains the global attractor. By symmetry this ball must be centered on the $z$-axis. Fixing its center to some $z$-coordinate $z_0$, we choose $\Phi=x^2+y^2+(z-z_0)^2$. The bound on $\Phi_\cA^+$ found by solving \eqref{eq: opt sos} give the squared radius of a sphere centered at $z_0$ which contains $\cA$, and we can solve \eqref{eq: opt sos} repeatedly while searching over $z_0$ to find the $z_0$ value that minimizes this radius. Using $V$ of degree 6, we find that the ball of radius 32.7044 centered at $z_0=29.9587$ contains the global attractor. This ball is shown in \cref{fig: sphere}, along with numerically integrated trajectories starting along the unstable manifold of the origin. Repeating the computations with $V$ of degree 8 gives a ball whose radius is barely smaller, differing only after the fifth digit.

\subsection{\label{sec: transient}Transient chaos}

To test the quality of bounds computed using \eqref{eq: opt sos} for a system with a chaotic saddle, we consider the so-called nine-mode model of Moehlis \emph{et al.}~\cite{Moehlis2004, Moehlis2005}. It is a modal truncation of the Navier--Stokes equations that captures several qualitative features in the regime of transition to turbulence. The system has a state vector $\mathbf a\in\R^9$ of mode amplitudes, three geometric parameters ($\alpha,\beta,\gamma$), and a parameter called the Reynolds number ($\Rey$). Each ODE takes the form
\begin{equation}
\ddt a_i = - \frac{d_i}{\Rey}(a_i-\delta_{i1}) + \sum_{j,k=1}^9N_{ijk} a_j a_k,
\label{eq: 9 mode}
\end{equation}
where $\delta_{ij}$ is the Kronecker delta, and the coefficients $d_i$ and $N_{ijk}$ depend on the geometric parameters. The system of nine equations is given explicitly in \Cref{app: 9 mode}.

For the present example we fix the parameters to $(\alpha,\beta,\gamma)=(1/2,\pi/2,1)$ and $\Rey=105$. At these parameters, trajectories that are numerically integrated from various initial conditions display moderately long chaotic transients before eventually tending to the fixed point at $\mathbf a=e_1$, where $e_1$ is the unit vector in the direction of the first coordinate \cite{Moehlis2004}. The fixed point is not globally attracting since there exist other invariant structures \cite{Moehlis2005, Lakshmi2020}, but it does appear to attract generic initial conditions. This behavior suggests that the system has a chaotic saddle---a complicated invariant set that is not an attractor but is a subset of the global attractor $\cA$. In such cases, a trajectory starting away from the global attractor is expected to have a short initial transient during which it approaches the saddle, followed by the much longer chaotic transient during which it shadows the chaotic saddle. In practice the chaotic phase of the trajectory is so close to the saddle that it lies in $\cA$ up to any reasonable numerical precision, in which case our bounds on $\Phi_\cA^+$ apply to this chaotic phase.

Quantities of physical interest in the nine-mode model include the perturbation energy $\mathcal E$ and total dissipation $\mathcal D$,
\begin{align}
\mathcal E &:= |\mathbf a-e_1|^2, &
\mathcal D &:= \frac{1}{Re}\sum_{i=1}^9 d_ia_i^2,
\end{align}
The simple steady state $\mathbf a=e_1$ minimizes $\mathcal E$ and maximizes $\mathcal D$. On the other hand, relatively large values of $\mathcal E$ and small values of $\mathcal D$ are found on complicated trajectories. Here we compute upper bounds on $\mathcal E$ and lower bounds on $\mathcal D$ over the global attractor by solving \eqref{eq: opt sos}. Before doing so, we can impose some structure on $V$, as described in \cref{sec: structure}. First, in order for $V$ and $\f\cdot\nabla V$ to both have highest-degree terms of even degree, which is necessary to satisfy the SOS constraints, we require the highest-degree terms of $V$ to take the form $|\mathbf a|^{2p}(a_1-a_9)^{2q}$. Second, the quantities $\mathcal E$ and $\mathcal D$ to be bounded are invariant under the same two symmetries as the ODEs, which are~\cite{Lakshmi2020}
\beq
\begin{array}{l}
(a_1,a_2,a_3,a_4,a_5,a_6,a_7,a_8,a_9)\mapsto
	(a_1,a_2,a_3,-a_4,-a_5,-a_6,-a_7,-a_8,a_9), \\
(a_1,a_2,a_3,a_4,a_5,a_6,a_7,a_8,a_9)\mapsto
	(a_1,-a_2,-a_3,a_4,a_5,-a_6,-a_7,-a_8,a_9).
\end{array}
\eeq
Thus we can impose these same symmetries on $V$.

\Cref{tab: 9 mode} reports the upper bounds on $\mathcal E$ and lower bounds on $\mathcal D$ that we have computed by solving \eqref{eq: opt sos} for $V$ of polynomial degrees up to 8. (To avoid small quantities, we report results using the scaled dissipation $\widehat{\mathcal D}:=Re\,\mathcal D$.) The bounds improve with each increase in degree and have not yet converged at degree 8. We do not report results of degree-10 computations, which were inaccurate due to poor numerical conditioning.

\begin{table}[tp]
\caption{\label{tab: 9 mode}Upper bounds on the maximum perturbation energy ($\mathcal E$) and lower bounds on the minimum scaled dissipation ($\widehat{\mathcal D}=Re\,\mathcal D$) over the global attractor of the nine-mode model. Bounds are computed by solving \eqref{eq: opt sos} with $V$ of polynomial degrees 2, 4, 6, and 8. The values of $\lambda$ giving the tabulated bounds are reported by \cref{tab: 9 mode lambda} in \Cref{app: results}. Also shown are the maximum $\mathcal E$ and minimum $\widehat{\mathcal D}$ found by numerical integration (see text).}
\begin{center}
\begin{tabular}{cccccc}
$\Phi$ & \multicolumn{4}{c}{Bounds} & Integration\\ \hline
& $\deg=2$ & $\deg=4$ & $\deg=6$ & $\deg=8$\ \\ \cline{2-5}
$\mathcal E$ & 2.146~\, & 1.454 & 1.160 & 1.112 & 1.039 \\
$\widehat{\mathcal D}$ & 0.0001 & 0.030 & 0.130 & 0.218 & 0.268
\end{tabular}
\end{center}
\end{table}

To judge the sharpness of our bounds, we sought trajectories on the chaotic saddle that attain large $\mathcal E$ or small $\mathcal D$. The dimension of state space makes this search much harder than for the Lorenz system in the previous subsection; some unstable periodic orbits of the nine-mode model have been computed \cite{Moehlis2005, Lakshmi2020} but not many, and their unstable manifolds have not been examined. However, a recently discovered periodic orbit appears to both maximize the infinite-time average of $\mathcal E$ and minimize the infinite-time average of $\mathcal D$ when $\Rey=105$ \cite{Lakshmi2020}. We focus our search on the unstable manifold of this orbit, called PO1 in \cite{Lakshmi2020}, which indeed gives relatively large maximum $\mathcal E$ and small  minimum $\mathcal D$. To explore the unstable manifold, we chose a point on PO1 and numerically integrated $2\cdot 10^4$ trajectories starting from small random perturbations of this point. The largest $\mathcal E$ and smallest $\widehat{\mathcal D}$ found among the resulting trajectories are reported in \cref{tab: 9 mode}. These values agree in their first digits with the corresponding degree-8 bounds, confirming that the bounds are fairly tight. We cannot say exactly how tight the bounds are because the extreme values found by numerical integration are unlikely to be true extrema over the global attractor. Indeed, the great difficulty of searching over the global attractor is what makes \emph{a priori} bounds valuable.

\subsection{\label{sec: multi}Multiple basins of attraction}

\begin{figure}[tp]
\begin{center}
\includegraphics[width=130pt,trim={26pt 0pt 58pt 18pt},clip]{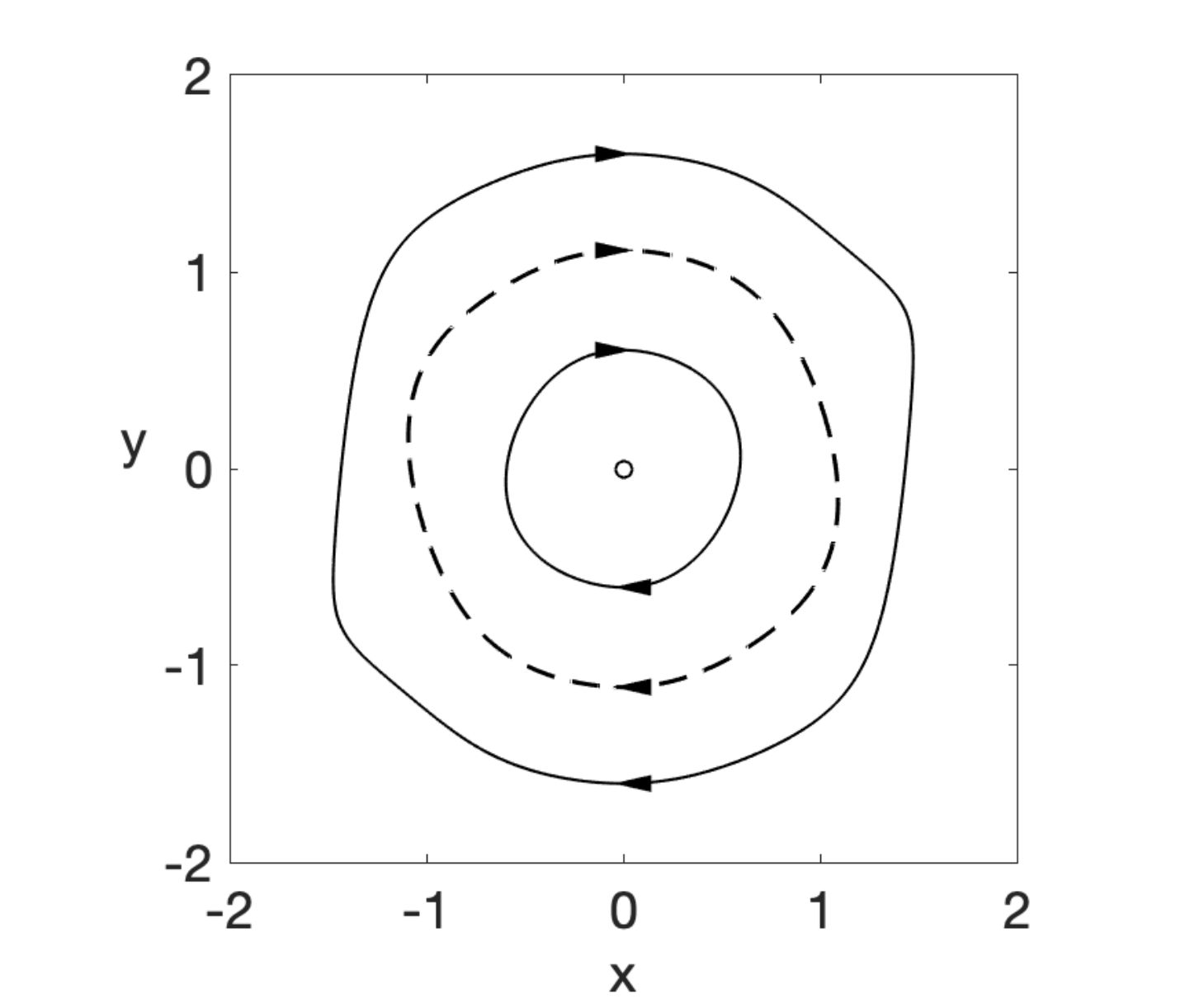}
\end{center}
\caption{\label{fig: multi}Phase portrait of \eqref{eq: multi} showing the attracting (\solidrule) and repelling (\dashedrule) limit cycles and the repelling equilibrium ($\circ$).}
\end{figure}

To illustrate how our non-global framework \cref{eq: opt X sos} applies to systems with multiple basins of attraction, we consider a two-dimensional example from \cite{Gasull2002} that has two locally attracting limit cycles. For this system, the components of the generic ODE \eqref{eq: ode} are
\beq
\x = (x,y), \quad 
\f = \big(y-x(x^2-2)(x^2-1)(x^2-1/4) ,\,-x\big).
\label{eq: multi}
\eeq
As shown by the phase portrait in \cref{fig: multi}, the basins of the two attracting limit cycles are separated by a repelling limit cycle, and there is a repelling equilibrium at the origin.

Suppose we are interested only in trajectories starting inside the middle limit cycle. All such trajectories eventually remain in the ball $X$ of radius 4/5 centered at the origin. (This $X$ contains the inner limit cycle and does not intersect the middle limit cycle.) The minimal set $\cA_X$ that attracts all bounded subsets of $X$ is composed of the inner limit cycle and the area within it. We seek upper bounds on the maxima of $\Phi=x^2$, $xy$, and $y^2$ over $\cA_X$. Sharpness of such bounds is easy to judge since each maximum is attained on the inner limit cycle and can be found by numerical integration. It would be impossible for the global formulation \eqref{eq: opt sos} to give a sharp bound on $\Phi_{\cA_X}^+$ in these cases since the global maxima $\Phi_\cA^+$ are attained on the outer limit cycle and are strictly larger than $\Phi_{\cA_X}^+$.

We computed upper bounds on the maxima of $\Phi=x^2$, $xy$, and $y^2$ by solving \eqref{eq: opt X sos} with $g(x,y)=(4/5)^2-x^2-y^2$. For $V$ of polynomial degrees up to 10, we performed computations with $\deg(s_i)=\deg(V)-2$ and $\deg(s_i)=\deg(V)$. Bounds are reported only for the latter choice, which improves bounds by enough to merit the additional computational cost. The ODE was rescaled by $(x,y)\mapsto1.6(x,y)$ so that all limit cycles lie in $[-1,1]^2$, which improves numerical conditioning, and then the computed bounds were multiplied by $1.6^2$ to recover the original scaling. The $\Phi$ to be bounded share the symmetry $(x,y)\mapsto(-x-y)$ of the ODE, so we can impose this same symmetry on $V$ \emph{a priori}, as described in \cref{sec: structure}.

\begin{figure}[tp]
\begin{center}
\includegraphics[width=220pt,trim={20pt 2pt 42pt 16pt},clip]{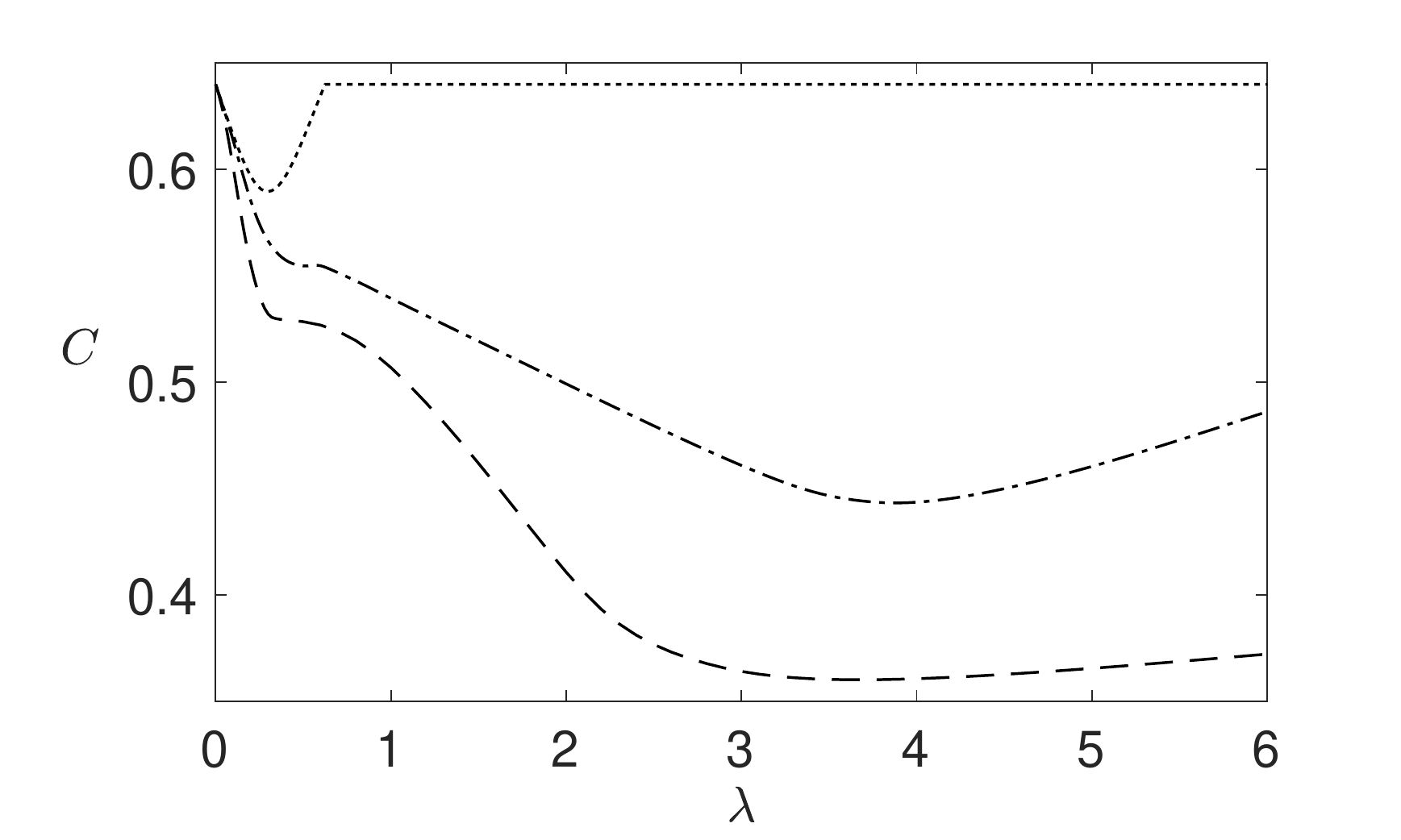}
\end{center}
\caption{\label{fig: lambda X}Upper bounds $C$ on the maximum of $x^2$ over the minimal set $\cA_X$ that attracts all subsets of the ball $X$ of radius 4/5 for the system \eqref{eq: multi}. The bounds are optima of the inner minimization in \eqref{eq: opt X sos} for various $\lambda$ and for $V$ of polynomial degrees 2 (\dottedrule), 4 (\dashdottedrule), and 6 (\longdashedrule).}
\end{figure}

The optimal bounds given by the inner minimization problem in \eqref{eq: opt X sos} display continuous but non-convex dependence on $\lambda$, unlike the convex dependence we observed for the global problem that is typified by \cref{fig: lambda}. As an example in the non-global case, \cref{fig: lambda X} shows how upper bounds on the maximum of $x^2$ over $\cA_X$ depend on $\lambda$. The fact that $X$ contains $\cA_X$ gives the trivial bound $\max_{\cA_X}x^2\le\max_X x^2 = (4/5)^2$. For $V$ of degrees 2, 4, and 6, bounds on $x^2$ are better than this trivial bound when $\lambda$ is less than approximately 0.624, 12.7, and 234, respectively.

\Cref{tab: multi} reports the optimal bounds we find by carrying out the non-convex but one-dimensional minimization over $\lambda$ in \eqref{eq: opt X sos}. Also shown are numerical approximations of each $\Phi_{\cA_X}^+$, found by numerically integrating a trajectory that tends to the inner limit cycle. Bounds improve as the degrees of $V$ and $s_i$ are raised, and at degree 10 they are sharp to at least 4 digits.

\begin{table}[tp]
\caption{\label{tab: multi}Upper bounds over the minimal set $\cA_X$ that attracts all subsets of the ball $X$ of radius 4/5 for the system \eqref{eq: multi}, computed by solving \eqref{eq: opt X sos} with $V$ of degrees up to 10. In all cases $\deg(s_i)=\deg(V)$.  Also shown is a numerical approximation of each maximum, $\Phi_{\cA_X}^+$, which occurs on the inner limit cycle shown in \cref{fig: multi}. The values of $\lambda$ giving the tabulated bounds are reported by \cref{tab: multi lambda} in \Cref{app: results}.}
\begin{center}
\begin{tabular}{ccccccc}
$\Phi$ & \multicolumn{5}{c}{Upper bounds on $\Phi_{\cA_X}^+$} 
	& $\Phi_{\cA_X}^+$ \\[2pt] \hline
& $\deg=2$ & $\deg=4$ & $\deg=6$ & $\deg=8$ & $\deg=10$ \\ \cline{2-6}
$x^2$ & 0.5896 	& 0.4432 & 0.3601 & 0.3555 & 0.3552 & 0.3552 \\
$xy$   & 0.3194 	& 0.2505 & 0.1957 & 0.1928 & 0.1926 & 0.1926 \\
$y^2$ & $(4/5)^2$ 	& 0.4608 & 0.3681 & 0.3634 & 0.3631 & 0.3631
\end{tabular}
\end{center}
\end{table}

\section{\label{sec: con}Conclusions}

We have presented a method for bounding extreme values of chosen quantities over global attractors---that is, over the minimal set $\cA$ that attracts all bounded sets of initial conditions. We also have presented a non-global version for bounding extreme values over $\cA_X$, the minimal set that attracts all bounded subsets of a specified set $X$. Our approach involves constructing Lyapunov functions by solving convex optimization problems. When all quantities are polynomial, a version of the framework can be implemented computationally by solving polynomial optimization problems subject to sum-of-squares constraints.

To illustrate our computational approach, we have reported bounds on various quantities for three ODE examples: the chaotic Lorenz system, a nine-mode model displaying transient chaos, and a system with multiple basins of attraction. Bounds over global attractors were computed for the first two examples, and bounds over a subset of the global attractor in a particular basin were computed for the third example. For every quantity bounded in the Lorenz system and the multiple-basin example, the computed bounds are very sharp and appear to become arbitrarily sharp as the polynomial degree of Lyapunov functions is raised. Most bounds in the prior literature on the Lorenz equations are not nearly as sharp. In the transient chaos example we cannot judge the sharpness of bounds well enough to say that they become arbitrarily sharp, nor is there any indication to the contrary. Preliminary computations (not reported here) for some other ODEs, including the Lorenz-84 model \cite{Lorenz1984} and truncations of the Kuramoto--Sivashinsky equation, also appear to give arbitrarily sharp bounds with increasing polynomial degree.

A fundamental theoretical question that remains open is: under what conditions does the present method produce sharp bounds? In particular, when are the maxima defining the left-hand quantities in \cref{eq: opt general,eq: opt X} equal to the minima of the right-hand problems? Results of this type have been proved for similar variational methods which give bounds on infinite-time averages \cite{Tobasco2018} or extreme events \cite{Fantuzzi2019}. Sharpness of the formulation \cref{eq: opt X} would imply sharpness of its computational relaxation \eqref{eq: opt X sos} in the polynomial case, as discussed at the end of \cref{sec: sos}, thereby ensuring that our successful examples are in fact typical. Another theoretical question of practical importance is: under what conditions do the inner minima in \eqref{eq: opt sos} or \eqref{eq: opt X sos} have convex dependence on $\lambda$, thereby simplifying the outer minimization problems?

The results reported here constitute yet another instance where methods based on polynomial optimization, when applicable, produce stronger results about dynamical systems than any other approach. Related methods have been similarly successful for tasks such as demonstrating stability \cite{Parrilo2000, Papachristodoulou2002, Anderson2015a}, bounding time averages \cite{Chernyshenko2014, Fantuzzi2016, Goluskin2018}, and estimating basins of attraction \cite{Tan2008, Chesi2011, Henrion2014, Anderson2015a, Valmorbida2017, Drummond2018}. Application to high-dimensional dynamical systems remains a practical challenge that calls for improving scalability, perhaps by replacing sum-of-squares constraints with stronger constraints that are more computationally tractable \cite{Waki2006, Ahmadi2017a, Ahmadi2017c, Ahmadi2019a, Zheng2019b}. Nevertheless, the further development of polynomial optimization methods for ordinary differential equations is sure to remain fruitful, as is the extension of such methods to partial differential equations.

\vspace{-10pt}
\paragraph*{Acknowledgements}
The author thanks Giovanni Fantuzzi for valuable discussions while visiting the 2018 Geophysical Fluids Dynamics program at the Woods Hole Oceanographic Institution. Helpful remarks by Sergei Chernyshenko and Charles Doering are appreciated also. Thanks to Mayur Lakshmi for sharing data on periodic orbits of the nine-mode model. This work was supported by the NSERC Discovery Grants Program through awards RGPIN-2018-04263, RGPAS-2018-522657, and DGECR-2018-00371.

\appendix

\section{\label{app: results}Optimal $\lambda$ values}

\Crefrange{tab: UB lambda}{tab: multi lambda} give values of $\lambda$ that yield the bounds reported in \crefrange{tab: UB}{tab: multi}, respectively. The first three tables report solutions of the global bounding formulation \eqref{eq: opt sos}, while the last reports solutions of the non-global formulation \eqref{eq: opt X sos}. We believe these $\lambda$ values are globally optimal, but we cannot be certain because the optimizations over $\lambda$ in \eqref{eq: opt sos} and \eqref{eq: opt X sos} are non-convex in general. Bounds are fairly insensitive to the value of $\lambda$, especially when the degree of $V$ is large. This can be seen in \cref{fig: lambda,fig: lambda X}. Thus, to compute bounds that are optimal up to the precision reported in \crefrange{tab: UB}{tab: multi}, less precision is needed in $\lambda$. This is reflected in the precision of the $\lambda$ values reported in \crefrange{tab: UB lambda}{tab: multi lambda}, which is low but still sufficient to give the optimal bounds.

\begin{table}[tp]
\caption{\label{tab: UB lambda}Values of $\lambda$ for which we find the minimum upper bounds in \eqref{eq: opt sos} when bounding various normalized monomials over the global attractor of the Lorenz equations using $V$ of degrees 2, 4, and 6. The bounds computed using these $\lambda$ appear in \cref{tab: UB}.}
\begin{center}
\begin{tabular}{cccc}
Monomial & \multicolumn{3}{c}{$\lambda$} \\ \hline
& $\deg=2$ & $\deg=4$ & $\deg=6$ \\ \cline{2-4}
$x$ & 0.566 & 0.374 & 0.375 \\
$y$ & 0.5 & 0.375 & 0.375 \\
$z$ & 0.5 & 0.378 & 0.375 \\
$xy$ & 0.554 & 0.378 & 0.375 \\
$xz$ & 0.597 & 0.374 & 0.375 \\
$yz$ & 0.5 & 0.375 & 0.375 \\
$x^2y$ && 0.374 & 0.375 \\
$x^2z$ && 0.379 & 0.375 \\
$xy^2$ && 0.374 & 0.375 \\
$xyz$ && 0.379 & 0.375 \\
$xz^2$ && 0.375 & 0.375 \\
$y^2z$ && 0.380 & 0.375 \\
$yz^2$ && 0.375 & 0.375
\end{tabular}
\end{center}
\end{table}

\begin{table}[tp]
\caption{\label{tab: LB lambda}Values of $\lambda$ for which we find the maximum lower bounds in \eqref{eq: opt sos} when bounding various normalized monomials over the global attractor of the Lorenz equations using $V$ of degrees 2, 4, 6, and 8. The bounds computed using these $\lambda$ appear in \cref{tab: LB}.}
\begin{center}
\begin{tabular}{ccccc}
Monomial & \multicolumn{4}{c}{$\lambda$} \\ \hline
& $\deg=2$ & $\deg=4$ & $\deg=6$ & $\deg=8$ \\ \cline{2-5}
$z$ & 0.375 \\
$xy$ & 0.554 & 0.408 & 0.375 & 0.375 \\
$x^2z$ && 0.388 & 0.296 & 0.321 \\
$xyz$ && 0.425 & 0.376 & 0.375 \\
$y^2z$ && 0.350 & 0.312 & 0.285 
\end{tabular}
\end{center}
\end{table}

\begin{table}[tp]
\caption{\label{tab: 9 mode lambda}Values of $\lambda$ for which we find the minimum upper bounds on $\mathcal E$ and maximum lower bounds on $\widehat{\mathcal D}$ for the nine-mode model of \cref{sec: transient}. Bounds are computed by solving \eqref{eq: opt sos} with $V$ of degrees 2, 4, 6, and 8. The bounds computed using these $\lambda$ appear in \cref{tab: 9 mode}.}
\begin{center}
\begin{tabular}{ccccc}
$\Phi$ & \multicolumn{4}{c}{$\lambda$} \\ \hline
& $\deg=2$ & $\deg=4$ & $\deg=6$ & $\deg=8$\ \\ \cline{2-5}
$\mathcal E$ & 32.8 & 43.1 & 48.3 & 47 \\
$\widehat{\mathcal D}$ & 27.6 & 29.7 & 16.2 & 21 \\
\end{tabular}
\end{center}
\end{table}

\begin{table}[tp]
\caption{\label{tab: multi lambda}Values of $\lambda$ for which we find the minimum upper bounds on $\Phi_{\cA_X}^+$ by solving \eqref{eq: opt X sos} for the system described in \cref{sec: multi}, with $V$ and $s_i$ of degrees up to 10. The bounds computed using these $\lambda$ appear in \cref{tab: multi}.}
\begin{center}
\begin{tabular}{cccccc}
$\Phi$ & \multicolumn{5}{c}{$\lambda$} \\ \hline
& $\deg=2$ & $\deg=4$ & $\deg=6$ & $\deg=8$ & $\deg=10$ \\ \cline{2-6}
$x^2$ & 0.297 & 3.88 & 3.67 & 3.42 & 3.38 \\
$xy$   & 0.052 & 3.84 & 3.71 & 3.45 & 3.37 \\
$y^2$ &	       & 3.72 & 3.59 & 3.44 & 3.37
\end{tabular}
\end{center}
\end{table}

\newpage

\section{\label{app: 9 mode}Nine-mode model}

We have used the nine-mode model of Moehlis \emph{et al.}~\cite{Moehlis2004} for the example of \cref{sec: transient}. To give the governing ODEs explicitly, recall that the system has four parameters ($\alpha,\beta,\gamma,\Rey$) and let
\begin{align}
\kappa_{\alpha \gamma} &= \sqrt{\alpha^2 + \gamma^2}, &
\qquad \kappa_{\beta \gamma} &= \sqrt{\beta^2 + \gamma^2}, &
\qquad \kappa_{\alpha \beta \gamma} &= \sqrt{\alpha^2 + \beta^2 + \gamma^2}.
\end{align}
The ODEs are~\cite{Moehlis2004}: 
%

\begin{align*}
\frac{{\rm d}a_1}{{\rm d}t} =& -\frac{\beta^2}{\Rey}(a_1-1) 
- \sqrt{\frac{3}{2}} \frac{\beta \gamma}{\kappa_{\alpha \beta \gamma}} a_6 a_8 
+  \sqrt{\frac{3}{2}}\frac{\beta \gamma}{\kappa_{\beta \gamma}} a_2 a_3, \\
\frac{{\rm d}a_2}{{\rm d}t} =& -\left( \frac{4 \beta^2}{3} + \gamma^2 \right) 
\frac{a_2}{\Rey} 
+ \frac{5 \sqrt{2}}{3 \sqrt{3}} \frac{\gamma^2}{\kappa_{\alpha \gamma}} a_4 a_6 
- \frac{\gamma^2}{\sqrt{6} \kappa_{\alpha \gamma}} a_5 a_7 \nonumber \\
& - \frac{\alpha \beta \gamma}{\sqrt{6} \kappa_{\alpha \gamma} 
\kappa_{\alpha \beta \gamma}} a_5 a_8
 - \sqrt{\frac{3}{2}} \frac{\beta \gamma}{\kappa_{\beta \gamma}} a_1 a_3 
- \sqrt{\frac{3}{2}} \frac{\beta \gamma}{\kappa_{\beta \gamma}} a_3 a_9, \\
\frac{{\rm d}a_3}{{\rm d}t} =& -\frac{\beta^2 + \gamma^2}{\Rey} a_3 
+ \frac{2}{\sqrt{6}} \frac{\alpha \beta \gamma}{\kappa_{\alpha \gamma} 
\kappa_{\beta \gamma}} (a_4 a_7 + a_5 a_6) \nonumber 
+ \frac{\beta^2 (3 \alpha^2 + \gamma^2)
- 3 \gamma^2 ( \alpha^2 + \gamma^2)}{\sqrt{6} \kappa_{\alpha \gamma} 
\kappa_{\beta \gamma} \kappa_{\alpha \beta \gamma}} a_4 a_8, \\
\frac{{\rm d}a_4}{{\rm d}t} =& -\frac{3 \alpha^2 + 4 \beta^2}{3 \Rey} a_4 
- \frac{\alpha}{\sqrt{6}} a_1 a_5 
- \frac{10}{3 \sqrt{6}} \frac{\alpha^2}{\kappa_{\alpha \gamma}} a_2 a_6 
\nonumber \\
& - \sqrt{\frac{3}{2}} \frac{\alpha \beta \gamma}{\kappa_{\alpha \gamma} 
\kappa_{\beta \gamma}} a_3 a_7 - \sqrt{\frac{3}{2}} 
\frac{\alpha^2 \beta^2}{\kappa_{\alpha \gamma} 
\kappa_{\beta \gamma} \kappa_{\alpha \beta \gamma}} a_3 a_8
- \frac{\alpha}{\sqrt{6}} a_5 a_9, \\
\frac{{\rm d}a_5}{{\rm d}t}
=& -\frac{\alpha^2 + \beta^2}{\Rey} a_5 + \frac{\alpha}{\sqrt{6}} a_1 a_4 
+ \frac{\alpha^2}{\sqrt{6} \kappa_{\alpha \gamma}} a_2 a_7 \nonumber \\
&- \frac{\alpha \beta \gamma}{\sqrt{6} \kappa_{\alpha \gamma} 
\kappa_{\alpha \beta \gamma}} a_2 a_8 + \frac{\alpha}{\sqrt{6}} a_4 a_9 
+ \frac{2}{\sqrt{6}} \frac{\alpha \beta \gamma}{\kappa_{\alpha \gamma} 
\kappa_{\beta \gamma}} a_3 a_6, \\
\frac{{\rm d}a_6}{{\rm d}t} =& -\frac{3 \alpha^2 + 4 \beta^2 + 3 \gamma^2}{3 \Rey} a_6 + 
\frac{\alpha}{\sqrt{6}} a_1 a_7  + 
\sqrt{\frac{3}{2}} \frac{\beta \gamma}{\kappa_{\alpha \beta \gamma}} a_1 a_8 
\nonumber \\ 
& + \frac{10}{3 \sqrt{6}} \frac{\alpha^2 - \gamma^2}{\kappa_{\alpha \gamma}} a_2 a_4 
- 2 \sqrt{\frac{2}{3}} \frac{\alpha \beta \gamma}{\kappa_{\alpha \gamma} 
\kappa_{\beta \gamma}} a_3 a_5 + \frac{\alpha}{\sqrt{6}} a_7 a_9 + \sqrt{\frac{3}{2}} \frac{\beta \gamma}{\kappa_{\alpha \beta \gamma}} a_8 a_9, \\
\frac{{\rm d}a_7}{{\rm d}t} =& -\frac{\alpha^2 + \beta^2 + \gamma^2}{\Rey} a_7 
- \frac{\alpha}{\sqrt{6}} (a_1 a_6 + a_6 a_9) + \frac{1}{\sqrt{6}} \frac{\gamma^2 - \alpha^2}{\kappa_{\alpha \gamma}} a_2 a_5 
+ \frac{1}{\sqrt{6}} \frac{\alpha \beta \gamma}{\kappa_{\alpha \gamma} 
\kappa_{\beta \gamma}} a_3 a_4, \\
\frac{{\rm d}a_8}{{\rm d}t} =& -\frac{\alpha^2 + \beta^2 + \gamma^2}{\Rey} a_8 
+ \frac{2}{\sqrt{6}} \frac{\alpha \beta \gamma}{\kappa_{\alpha \gamma} 
\kappa_{\alpha \beta \gamma}} a_2 a_5 + \frac{\gamma^2 (3 \alpha^2 - \beta^2 + 3 \gamma^2)}{\sqrt{6} 
\kappa_{\alpha \gamma} \kappa_{\beta \gamma} \kappa_{\alpha \beta \gamma}} 
a_3 a_4, \\
\frac{{\rm d}a_9}{{\rm d}t} =& -\frac{9 \beta^2}{\Rey} a_9 + \sqrt{\frac{3}{2}} 
\frac{\beta \gamma}{\kappa_{\beta \gamma}} a_2 a_3 
- \sqrt{\frac{3}{2}} \frac{\beta \gamma}{\kappa_{\alpha \beta \gamma}} 
a_6 a_8.
\end{align*}

\bibliography{Global_attractor_bounds.bbl}

\end{document}